\magnification=\magstep1


\outer\def\give#1. {\medbreak
             \noindent{\sl#1. }}                     
\outer\def\section #1\par{\bigbreak\centerline{\S
     {\bf#1}}\nobreak\smallskip\noindent}
\def\({\left(}
\def\){\right)}
\def\contract{{\,{\vrule height5pt width0.4pt depth 0pt}
{\vrule height0.4pt width6pt depth 0pt}\,}}

\def\sqr#1#2{{\vcenter{\hrule height.#2pt              
     \hbox{\vrule width.#2pt height#1pt\kern#1pt
     \vrule width.#2pt}
     \hrule height.#2pt}}}
\def\square{\mathchoice\sqr{5.5}4\sqr{5.0}4\sqr{4.8}3\sqr{4.8}3}
\def\qed{\hskip4pt plus1fill\ $\square$\par\medbreak}




\def\cD{{\cal D}}

\def\cG{{\cal G}}

\def\cR{{\cal R}}

 
\def\CC{{\bf C}}

\def\C{{\bf C}}
\def\cx#1{{\CC}^{#1}}     


\def\bar{\overline}              
 
 


 
\def\C{{\bf C}}
\def\cx#1{{\CC}^{#1}}     

\centerline{\bf Distribution of Periodic Points}
\centerline{\bf of Polynomial Diffeomorphisms of $\cx2$}
\bigskip
\centerline{E. Bedford, M. Lyubich and J. Smillie}
\bigskip
 
\section 1. Introduction.
 
This paper deals with the dynamics of a simple family of  holomorphic
diffeomorphisms of $\cx 2$: the polynomial automorphisms.
 This family of maps has been studied by a number of authors. We
refer to [BLS] for a general introduction to this class of
dynamical systems.  An interesting object from the point of view of
potential theory is the
equilibrium measure $\mu$ of the set $K$ of points with bounded orbits.
In [BLS] $\mu$ is also characterized dynamically as the unique measure
of maximal entropy. Thus $\mu$ is also an equilibrium measure from
the point of view of the thermodynamical formalism.
  In the present paper we give another dynamical interpretation of
$\mu$ as the limit distribution of the periodic points of $f$.
 
Fix a polynomial automorphism $f$. A point $p$ is periodic if
$f^np=p$ for some $n>0$, and the smallest positive $n$ for which
this equation holds is the period of $p$.  We let ${\rm Fix}_n$
denote the set of fixed points of $f^n$, and  ${\rm Per}_n$ denote
the set of points of period exactly $n$.  Thus ${\rm Fix}_n
=\bigcup{\rm Per}_k$,  where the union is taken over all $k$ which
divide $n$.   The map $f$ has a dynamical degree $d$ which we assume to be
larger than 1. By [FM] $f^n$ has exactly $d^n$ fixed points counted with
multiplicities. Since the multiplicity of each fixed point of $f^n$ is positive
we conclude that $\#{\rm Fix}_n\le d^n$.
 
A periodic point $p$ is called a saddle point if the eigenvalues
$\lambda^s(p)$ and $\lambda^u(p)$ of $Df^n(p)$ satisfy
$|\lambda^s(p)|<1< |\lambda^u(p)|$.
 We let ${\rm SPer}_n$ denote the saddle points with period
exactly $n$, so that $${\rm SPer}_n\subset {\rm Per}_n\subset {\rm
Fix}_n\eqno(1)$$ $$\#{\rm SPer}_n\le \#{\rm Per}_n\le \#{\rm
Fix}_n\le d^n.\eqno(2)$$
 
\proclaim Theorem 1.   If $P_n$ denotes any of the three sets in
(1), then $$\lim_{n\to\infty}{1\over {d^n}}\sum_{a\in
P_n}\delta_a=\mu.$$
 
There seems to be a general tendency in many classes of dynamical systems for
periodic points to be equidistributed with respect to the measure of maximal
entropy. Indeed, this is true for subshifts of finite type and their smooth
counterparts, axiom A diffeomorphisms (Bowen [B1]).  This is also
true for  rational endomorphisms of the Riemann sphere [L].  In
the special  case of polynomial maps of the complex plane this
result can also be derived by Brolin's methods [Br] (see Sibony [Si] and Tortrat
[T].) For polynomial diffeomorphisms of $\C^2$ the equidistribution property was
conjectured by N. Sibony.  In [BS1], this was proven with the additional
hypothesis $f$ is hyperbolic (but without the a priori assumption that the
periodic points are dense in the Julia set).
 
Brolin uses potential theory to analyze polynomials of one
variable, where the harmonic measure is the unique measure
of maximal entropy. In the case of (non-polynomial) rational maps there is a
unique measure of maximal entropy but it is not in general related
to the harmonic measure. Lyubich's proof [L] in the case of rational
maps uses the balanced property of the unique measure of maximal
entropy together with a Shadowing Lemma. For polynomial
automorphisms of $\C^2$ the analogue of the balanced property is
the product structure as described in [BLS]. In this paper we use
the product structure and the Shadowing Lemma.  This argument is
of quite a general nature and we expect it to be useful outside
the holomorphic setting.  From our point of view the product
structure of the measure of maximal entropy is the
underlying reason for the equidistribution property in the known
cases.
 
According to Newhouse there exist polynomial automorphisms with
infinitely many sink orbits.  The next result shows that the
majority of orbits are saddle orbits.

\proclaim Corollary 1.  Most periodic points are saddle points in
the sense that $$\lim_{n\to\infty}{ 1 \over {d^n}}{\#{\rm SPer}_n}
=1.$$
 
A weaker asymptotic formula $$\limsup_{n\to\infty}\
(\#{\rm SPer}_n)^{1/n}=d $$ follows from a theorem of Katok [K]
and the entropy formula $h(f)= \log d$ ([FM] and [S]).
 
The following corollary answers a question of [FM]:
 
\proclaim Corollary 2. $f$ has points of all but finitely many
periods.

Given a point $x$, let us define {\it the Lyapunov exponent} at
$x$ as $$\chi(x)=\lim_{n\to +\infty}{1\over n}\log\|Df^n(x)\|,$$
provided the limit exists. For example, if $p$ is a saddle point
of period $n$, then $\chi(p)={1\over n}\log\lambda^u(p)$. For any
ergodic invariant measure $\nu$ the function $\chi(x)$ is constant
$\nu$ almost everywhere. This common value is called the Lyapunov
exponent of $f$ with respect to $\nu$. The harmonic measure
$\mu$ is ergodic, and we denote
by $\Lambda$
 the Lyapunov exponent of $f$ with respect to $\mu$.   We have an alternate
 description of $\Lambda$ as:
$$\Lambda=\lim_{n\to +\infty}{1\over n}\int\log||Df^n||\mu.$$
 
The following result allows us to compute $\Lambda$ by averaging
the Lyapunov exponents of the saddle points.

\proclaim Theorem 2. For any polynomial automorphism $f$ we have:
$$\Lambda=\lim_{n\to\infty}{1\over n{d^n}} \sum_{p\in {\rm
SPer}_n}\chi(p).$$
 
The statement is equally true with the set ${\rm SPer}_n$ replaced
by ${\rm Fix}_n$.  We note (from [BS3]) that
$\Lambda\ge\log d$, so the previous theorem gives a lower bound
for the average exponent for periodic points.
 
For quadratic maps in one dimension the Lyapunov exponent with
respect to the harmonic measure is closely related to the Green
function of the Mandelbrot set. In particular the dynamical
behavior of $f$ is reflected in the behavior of this function. It
is interesting to investigate the relation between the behavior of
$\Lambda=\Lambda(f)$ and the dynamics of $f$ for automorphisms of
$\cx2$.  Consider a holomorphic parameterization $c\mapsto f_c$.
It is shown in [BS3] that $\Lambda(f_c)$ is subharmonic for all
values of $c$. The following result shows that non-harmonicity of
$\Lambda(f_c)$ corresponds to the creation
 of sinks (in the dissipative case) or Siegel balls (in the
volume-preserving case) for nearby maps.  (A Siegel ball is the analogue for
volume-preserving maps of a Siegel disk, see [BS2], [FS].)
 
\proclaim Theorem 3. Consider a family of $\{f_c\}$ depending
holomorpically on a parameter $c$ in the disk. Assume that
$c\mapsto\Lambda(f_c)$ is not harmonic at $c=c_0$. If the maps
$f_c$ are dissipative, then there is a sequence $c_i\to c_0$
defined for $i\ge N$ so that $f_{c_i}$ has a sink of period $i$.
If the maps $f_c$ are volume preserving, then there is a sequence
$c_i$ defined for $i\ge N$ so that $f_{c_i}$ has a Siegel ball of period $i$.
 
It is shown in [BS3] that $c\mapsto\Lambda(f_c)$ is a harmonic
function for the values of the parameter for which $f_c$ is
hyperbolic. In this context hyperbolicity implies structural
stability which implies that the topological conjugacy class of
the map is locally constant. The proof of Theorem 3 shows that
structural stability directly implies the harmonicity of $\Lambda$:
 
\proclaim Theorem 4. If the maps $f_c$ are topologically conjugate
to one another in a neighborhood of $c=c_0$, then the function
$c\mapsto\Lambda(f_c)$ is harmonic at $c_0$.

\bigskip \section 2. Lyapunov charts, Pesin boxes and the
Shadowing Lemma.
 
\noindent{\bf Topological bidisks.} Let us take a pair $\cD^s$ and
$\cD^u$ of standard disks $\{z: |z|\leq 1\}$, and consider the
standard bidisk $\cD=\cD^u\times \cD^s$.
 We refer to the sets $\cD^s(b)=\cD^s\times\{b\}$ as  {\it stable
cross sections} and $\cD^u(a)=\{a\}\times\cD^u$ as  {\it unstable
cross sections}. Also, the boundary of $\cD=\cD^s\times\cD^u$ will
be partitioned into  $\partial\cD^s\cup\partial\cD^u$, where we set
$\partial^s\cD:=\cD^s\times\partial\cD^u$ and
 $\partial^u\cD:=\partial\cD^s\times\cD^u$.
 
Let us define a   {\it topological bidisk} $B$ as a compact set
homeomorphic to $\cD$ together with a homeomorphism
$h:B\rightarrow \cD^2$. Then $B$ inherits the structure of the
stable/unstable ($s/u$)
 cross sections $B^s(x)$ and $B^u(x)$, and the partition of the
boundary into $\partial^s B$ and $\partial^u B$, which is induced
by this homeomorphism. We note that the stable/unstable cross
sections are not going to be stable or unstable manifolds, but
rather an approximation to them. If $h$ is affine we can talk
about the {\it affine bidisk}.
 
Let us say that a two dimensional manifold  $\Gamma$ {\it $u$-overflows}
a bidisk $(B,h)$ if $h(\Gamma\cap B)\subset \cD$ is the graph of a
funcion  $\phi: \cD^u\rightarrow\cD^s$. If additionally
$\Gamma\subset B$ we also say that $\Gamma$ is $u$-{\it inscribed}
in $B$. By the slope of $\Gamma$ in $B$ we mean
$\max\|D\phi(z)\|$.  Similarly we can define the dual concept
related to the stable direction.

We will say that a topological bidisk $B_1$ {\it $u$-correctly
intersects a topological bidisk $B_2$} (or ``the pair $(B_1,B_2)$
intersects $u$-correctly") if every unstable cross-section
$B_1^u(x)$ $u$-overflows $B_2$ and every stable cross section
$B_2^s(x)$ $s$-overflows $B_1$.
 If the pair $(B_1,B_2)$ intersects $u$-correctly,  then the
intersection $B_1\cap B_2$ becomes a topological bidisk, if we
give it  the unstable cross sections from $B_1$ and the stable
cross sections from $B_2$. (In other words, the straightening
homeomeorphism is $h: x\mapsto (h_1^s(x), h_2^u(x))$ where
$h_i=(h_i^s,h_i^u): B_i\rightarrow \cD$ are straightening
homeomorphisms for $B_i$). In the case when $B_1\subset B_2$ the
$u$-overflowing property is equivalent to $\partial^s B_1\subset
\partial^s B_2$. Then we also say that $B_1$ is $u$-{\it
inscribed} in $B_2$.
 
 Of course, we have the
corresponding dual $s$-concepts.
 
\medskip \noindent{\bf Adapted Finsler metric and Lyapunov charts.}
We will use the exposition of Pugh and Shub [PS] as the reference
to the Pesin theory. Below we will adapt our presentation of the
theory to our specific goals.  Let us consider a holomorphic
diffeomorphism $f: \C^2\rightarrow \C^2$. Let  $\mu$ be  an
invariant, ergodic, hyperbolic measure  (the latter means that it
has non-zero characteristic exponents,  $\chi^s<0<\chi^u$).  Let
$\cR$ denote the set of  Oseledets regular points for $f$.
 For $x\in\cR$ there exists an invariant splitting of the tangent
space $T_x\cx 2=E^s_x\oplus E^u_x$ into a contracting direction
$E^s_x$ and an expanding direction $E^u_x$ which depend measurably
on $x$.
 Further, for any $0<\chi<\min(|\chi^s|,\chi^u)$ there is a
measurable function $C(x)=C_{\chi}(x)>0$ such that
$$\left|Df^{n}|_{E^s_x}\right|\le C(x)e^{-n\chi} \quad {\rm
and}\quad\left|Df^{-n}|_{E^u_x}\right|\le C(x)e^{-n\chi}\eqno (2)$$
for $n=1,2,3,\dots$ Set $Q_C=\{x: \; C(x)\leq C\}. $ These sets
exhaust a set of full $\mu$ measure as $C\to\infty$. Let us  fix a
big $C$ and refer to $Q_C$ as $Q$.
 
A key construction of the Pesin theory is a measurable change of
the metric  which makes $f$ uniformly  hyperbolic. Namely, for
$\theta\in (0,\chi)$ and $x\in \cR$,  set $$|v|^*=\sum_{n\geq 0}
e^{n\theta} |Df^{-n} v|\quad {\rm for}\quad v\in E^u_x,
  \eqno(3^u)$$ $$|v|^*=\sum_{n\geq 0} e^{n\theta} |Df^n v|\quad
{\rm for}\quad v\in E^s_x.  \eqno(3^s)$$ For arbitrary tangent
vector $v\in T_x$ set  $|v|^*=\max(|v^s|^*,|v^u|^*)$ where $v^s$
and $v^u$ are its stable and unstable components. This metric is
called {\it adapted Finsler} or {\it Lyapunov}. Observe that
$${1\over 2}|v|\leq |v|^*\leq B(x)|v| \eqno(4)$$ with a measurable
function $B(x)=C(x)/(1-e^{\theta-\chi})$. Hence the adapted
Finsler metric is equivalent to the Euclidean metric on the set
$Q$.
 
Let $r(x)>0$ be a measurable function. Then taking advantage of
natural identification of $\cx2$ with $T_x(\cx2)$, we can consider
a family of affine bidisks $$L(x)=\{x+v : |v|^*<r(x)\}. \eqno (5)$$
By the {\it inner size} of $L(x)$ we mean min$\{|v^s|, |v^u| :
x+v\in L(x)\}$. The Pesin theory provides us with a choice of the
``size-function" $r(x)$  with the following properties.
 
\smallskip \noindent (L1) The inner size of Lyapunov charts is bounded
away from 0 on
  the set $Q$.
Indeed, by (4) on this set
 the adapted Finsler metric is equivalent to the Euclidean metric.
  But then  $r(x)$ stays away from 0 as one can see from its
explicit
   definition in [PS], p.13.
 
For a complex one-dimensional manifold $\Gamma\subset
L(x)$,
 let us define its {\it cut-off iterate}, $f_x\Gamma$, as
$f\Gamma\cap L(fx)$ .
 
\smallskip
\noindent (L2) There is a measurable function $\kappa(x)>0$ with
the
  following property.  Denote by $\cG^u_x$ the family of complex
one-dimensional
   manifolds which are $u$-inscribed into $L(x)$ and have the slope
   less than $\kappa(x)$. Then the operation of cut-off iteration, $f_x$,
transforms $\cG^u_x$ into
   $\cG^u_{fx}$. In particular,
  the topological bidisk $f L(x)$ correctly intersects $L(x)$
  (where the bidisk structure on $f L(x)$ comes from $L(x)$).
 
\smallskip This allows us to repeat the cut-off process for
iterates of $f$.
 Let $f^n_x$ denotes the composition of the cut-off
  iterates at $x, fx,\dots, f^{n-1}x$.
   The hyperbolicity of $\mu$ implies:
 
 \smallskip
\noindent (L3) For any two manifolds $\Gamma_i\in \cG^u_x$,
  their cut-off iterates are getting close exponentially fast:
 $$C^1{\rm -dist} (f_x^n\Gamma_1, f^n_x\Gamma_2)\leq A(x)
e^{-\theta n} $$
  with a measurable  $A(x)$. Moreover, $A(x)$ is bounded on $Q$.
 
\smallskip
\noindent (L4) Let $p=f^{-n}x$. Then for any  manifold
  $\Gamma\in \cG^u_p$, and for any two points $y,z\in
f^n_p Gamma$,
  $${\rm dist} (f_p^{-n} y, f^{-n}_p z)\leq A(x) {\rm dist}(y,z)
e^{-\theta n}$$
   with $A(x)$ as above.

 \smallskip
(In the stable direction we of course have the corresponding
properties with
 the reversed time.)
 Such a family of affine bidisks $L(x)$  will be called   a family
 of {\it Lyapunov charts}.
 
 
\medskip\noindent{\bf Stable and unstable manifolds.}  Let us
consider the push-forward of a  Lyapunov chart and trim it down by
intersecting with the image chart.    Since by (L2)  the pair
$(fL(x), L(fx))$ intersects $u$-correctly,  the set
 $L_1^u(fx)= L(fx)\cap fL(x)$ a topological  bidisk $u$-inscribed
into $L(fx)$. Similarly,
 $L^s_1(x)= L(x)\cap f^{-1} L(fx)$ is a topological bidisk
$s$-inscribed into $L(x)$, and  $f: L^s_1(x)\rightarrow L_1^u(fx)$
is a bidisk diffeomorphism.

Property (L2) allows us continue inductively by noting that the
pair  $(L(f^nx),fL^u_{n-1})$ intersects $u$-correctly, and we can
consider topological bidisks  $$\eqalign{L^u_n(f^nx)=L(f^nx)\cap
fL^u_{n-1}(f^{n-1}(x)=\cr
  f^n \{ x : x\in L(x), fx\in L(fx),...,f^n x\in L(f^n x)\},\cr}
\eqno (6a)$$ and $$ L^s_n(x)=  f^{-n} L^u_n(f^nx)=
   \{ x: x\in L(x),  fx\in L(fx),..., f^n x\in L(f^n x)\}. \eqno
(6b)$$ Observe that $L^s_n(x)$ can also be defined as the
connected component of $L(x)\cap f^{-n} L(f^n x)$ containing $x$.
 
Furthemore, $ L^s_n(x)$ is $s$-inscribed into $L(x)$, and by the
dual property to (L3)  the sizes of $L^s_n(x)$ in the stable
direction  shrink down to zero. Hence the intersection
$$W^s_{loc}(x)=\bigcap_{n\geq 0} L_n^s(x)=\{x  : f^nx\in
L(f^nx),\; n=0,1,...\}   \eqno(7)$$ is a manifold of the family
$\cG^u_x$ $s$-inscribed into $L(x)$.  This is  exactly  the local
stable manifold through $x$.
 By the dual to (L4), for any $y,z\in W^s(x)$ $${\rm dist}(f^ny,
f^nz)\leq A(x)e^{-\theta n}$$ with a measurable function $A(x)$.
 
Reversing the time we obtain local unstable manifolds as well.
 
 \medskip\noindent {\bf Pesin boxes.}
 Recall that $Q=Q_C\subset\cR$ denotes the set on which $C(x)\leq
C$ in (2).
 This set is compact, and it may be shown that $Q\ni x\mapsto
E^{s/u}_x$  is continuous.  Further, $W^s(x)$ and $W^u(x)$
intersect transverally at $x$, and $W^{s/u}_{loc}(x)$ varies
continuously with $x\in Q$.
 Thus if $x,y\in Q$ are sufficiently close to each other,  then
there is a unique point of intersection  $[x,y]:=W^s_{loc}(x)\cap
W^u_{loc}(y)$.  It  follows that $[\,,] : Q\times Q\to \cx 2$ is
defined and continuous near the diagonal. Given a small
$\alpha>0$, let us set $$\bar Q=\bar Q_{C,\alpha}=\{[x,y]:\;
x,y\in Q,\; {\rm dist}(x,y)<\alpha\}.$$ This is a closed subset of
$\cR$ on which $C(x) < C/2$ in (2) (provided $\alpha$ is
sufficiently small), and which is locally closed with respect to
$[\,,]$-operation. By (L1) the inner size of the Lyapunov charts
$L(x)$ stays away from 0 for
 $x\in \bar Q$.  Let $\rho>0$ be a lower bound for this size.
 
 We say that a set $P$ has {\it product structure} if $P$ is
closed under $[\,,]$.  If in addition $P$ is (topologically)
closed and has positive measure, then it will be called a {\it
Pesin box}. Let us cover $\bar Q$ with countably many closed
subsets $X_i$ of positive  measure with diameter $<\alpha$. Taking
the $[\,,]$-closure of $X_i$, we obtain a covering of $\bar Q$
with a countably many Pesin boxes $P_i\subset \bar Q$. Since the
$\mu$ almost whole space can be exhausted by the sets $Q_C$, we
conclude that it can be covered by countably many Pesin boxes of
arbitrarily small diameter. Our next goal is to make these boxes
disjoint (at expense of losing a set of small measure).
 
\proclaim Lemma 1.  For any $\epsilon>0$ and $\eta>0$ there is a
finite family of disjoint Pesin boxes $P_i$ such that ${\rm
diam}(P_i)<\eta$ and  $\mu(\bigcup P_i)>1-\epsilon$.
 
\give Proof. As  we can cover a set of full measure by countably
many Pesin boxes, we can cover a set $X$ of measure $1-\epsilon$
by finitely many boxes.  Now let us make them disjoint by the same
procedure as used by Bowen
 (cf.\ [B2], Lemma 3.13).
 Namely, if two boxes $P_1$ and $P_2$ intersect, subdivide each of
them into four sets $P^j_i\subset P_i$
 with product structure in the following way: $$P^1_1=\{x\in P_1:
W^s(x)\cap P_2\ne\emptyset,W^u(x)\cap P_2\ne\emptyset\},$$
$$P^2_1=\{x\in P_1: W^s(x)\cap P_2=\emptyset,W^u(x)\cap
P_2\ne\emptyset\},$$ $$P^3_1=\{x\in P_1: W^s(x)\cap
P_2\ne\emptyset,W^u(x)\cap P_2=\emptyset\},$$ $$P^4_1=\{x\in P_1:
W^s(x)\cap P_2=\emptyset,W^u(x)\cap P_2=\emptyset\}.$$
 
Repeating this procedure, we can cover the set $X$ by a finite
number of disjoint sets $Q_i, i=1,\dots,N, $ with product
structure. However, $Q_i$ need not be closed.  To restore this
property, let us consider  closed sets $K_j\subset Q_j$ such that
 $\mu(Q_j-K_j)<\epsilon/N$. Completing these sets with respect to
the product structure, we obtain a suitable family of Pesin
boxes. \hbox{\qquad\qquad}\qed
 
\medskip\noindent{\bf A common chart.} Our goal is to have a
common Lyapunov chart for all points returning to a Pesin box.
 
Observe that the family of Lyapunov charts can be reduced in size
by a factor $t\leq 1$.
 Then the manifolds of the
family $\cG^u_x$ (correspondingly $\cG^s_x$) become almost
parallel within the scaled charts. In particular,
 the stable cross sections of the bidisks $L^s_n(x)$, as well as
the local stable manifolds truncated  by the scaled family of
charts  become    almost parallel.
 
Let $\rho$ be the lower bound of the size of the
Lyapunov charts $L(x),\; x\in Q$.  Let us take a Pesin box $P$ of size
$\eta<\rho/8$.  Let $a\in P$,  and let $B^s(a,r)\subset E^s_a$,
$B^u(a,r)\subset E^u_a$ denote  the Euclidean disks of radius $r$
centered at $a$ in the corresponding subspaces.
 Let us consider an affine bidisk $B=B^s(a,\rho/2)\times
B^u(a,\rho/2)$. Since the stable/unstable directions through the
points $x\in P$ are almost parallel, we have the inclusions:
$$P\subset B\subset \bigcap_{x\in P} L(x).$$
 
For $x\in P\cap f^{-n} P$ let $T=B_{n,x}^s$ be the component of
$B\cap f^{-n} B$ containing $x$.  Let also $y=f^nx$, $R=f^n
T=B^u_{-n,y}$.
 
\proclaim Lemma 2. The set $T$/(respectively $R$) is a topological
bidisk which  is $s/u$-correctly inscribed in $B$. The $s/u$-cross
sections of $T/R$ belong to the families $\cG^{s/u}$
respectively.  Moreover, $T$  $s$-correctly intersects $R$.
 
\give Proof.  For $z\in T$ let $T^u(z)=T\cap(z+E^u_a)$, and
similarly for $\zeta\in R$ let $R^s(\zeta)=R\cap (\zeta+E^s_a)$.
The first claim will follow from the following dual statements:
 
\smallskip \item {(i)} $f_x^n T^u(z)$ is a topological disk
$u$-correctly inscribed into
            $B$;
 
\smallskip \item {(ii)} $f_y^{-n} R^s(\zeta)$ is a topological disk
which is $s$-correctly inscribed into $B$.
 
\noindent Let us prove (i). Let us consider the Lyapunov chart
$L(x)\supset B$.   Then the disk $K=L(x)\cap (z+E^u_a)$ belongs to
the family  of graphs ${\cG}^u_x$. Hence its cut-off iterate
$f_x^n K\in\cG^u(y)$ overflows $L(y)$. Since it is almost parallel
to $E^u_a$,  $f_x^n T^u(z)=f_x^n K\cap B$ is $u$-correctly inscribed
into $B$. This proves the first two claims.
 
The last one now follows from the transversality of the families
$\cG^u$ and $\cG^s$. \qed

 
Set $W^{u/s}_B(x)=W^{u/s}_{loc}(x)\cap B$.
 
\proclaim Lemma 3. Under the above circumstances we have
$$W^s_B(x)\subset T   \eqno (8)$$ and  $$W^u_{loc}(x)\cap T\subset
f^{-n} W^u_{loc}(y). \eqno (9)$$
 
\give Proof.   
 Since $f^n W^s_B(x)$ has an exponentially small size, it is
contained in $B$.
  So $W^s_B(x)\subset B\cap f^{-n}B$.  Since $W^s_B(x)$ is
connected, (8) follows. \qed
 
 In order to get (9),
observe that $T\subset L^s_n(x)$. But we know by the construction
of the stable/unstable manifolds that $ W^u_{loc}(x)\cap
L^s_n(x)=f^{-n} W^u_{loc}(y)$.   \qed

Observe that by definition the sets $B^u_{n,x}$ are either
disjoint or coincide. So, for each $n$ we can consider the
following equivalence  relation on the set $P\cap f^{-n} P$ of
returning points:
  $x\sim y$ if $B^u_{n,x}= B^u_{n,y}$.
 
\proclaim Lemma 4. The equivalence classes have a product
structure.
  Moreover, if $x,z\in P\cap f^{-n}P$ are equivalent then
  $f^n[x,z]=[f^nx,f^nz]$.
 
\give Proof. Denote $q=[x,z]$. Clearly, $f^n q\in W^s_{loc}(f^nx)$.
 Further, by (8), $q\in T$. Hence by (9) $f^nq\in W^u_{loc}(f^nz)$,
  and the property $f^n[x,z]=[f^nx,f^nz]$ follows. Moreover, it
follows
   that  $q\in T\cap P\cap f^{-n}P$, so that it is equivalent to
$x,z$.
                                                      \qed

\medskip \proclaim The Shadowing Lemma. For any $x\in P\cap
f^{-n}P$, $y=f^nx$,
  there is a unique
  saddle point $\alpha\in B^s_{n,x}\cap B^u_{-n,y}$ of period $n$.
  Moreover, ${\rm dist}(\alpha, P)\leq Ce^{-n\theta}$.
 
\give Proof. Let us consider the family $\cG^u$ of manifolds
$u$-inscribed into
  $B$ with the slope $\leq\kappa\leq\inf_{x\in P}\kappa(x)$
  (with respect to the decomposition $E^u_a\oplus E^s_a$) .
  Take a returning point
   $x\in P\cap f^{-n}P$, and let $y=f^nx$. Then  we can consider
   the cut-off iterate
    $\Phi_u: \Gamma\mapsto f^n_x \Gamma\cap B$. Since its image
$\Phi^u\Gamma$
   is close to the unstable manifold $W^u_B(y)$, it is almost
parallel
    to $E^u_a$. It follows that for sufficiently big $n$,
    $\Phi_u$ maps $\cG^u$ into itself. Moreover, this
transformation is
    contracting according to (L3).

   Hence there is a unique $\Phi_u$-invariant manifold
$G^u\in\cG^u$. This manifold is $u$-inscribed into $B^u_{-n,y}$,
and can be characterized as the set of all points non-escaping
$B^u_{-n,y}$ under backward  iterates of $f^n$. Similarly,  we can
define a transformation $\Phi_s$ corresponding to the return of
$y$ back to $P$ under $f^{-n}$, and find a unique
$\Phi_s$-invariant manifold $\Gamma^s\in\cG^s$, $\Gamma^s\subset
B^s_{n,x}$. This manifolds intersect transversally at a unique
point
 $\alpha$ which is the desired periodic point.
 
There are no other periodic ponts of period $n$ in $B^u_{x,n}$.
Indeed, all points escape $B^u_{x,n}$ either under forward or
backward  iterates of $f^n$. Finally since the bidisks $B^s_{n,x}$
and $B^u_{-n,y}$ are exponentially thin, the point $\alpha$ is
exponentially close to  $[x,y]\in P$.
 \qed
 
\bigskip \section 3. Proofs of the theorems.
 
According to the Shadowing Lemma,   to each returning point $x\in
P\cap f^{-n}P$ we can assign a periodic point
$\alpha=\alpha(x)\in B^s_{n,x}\cap B^u_{n,x}$ of period $n$.
Given such an $\alpha$, set $$T(\alpha)=\{x\in P\cap
f^{-n}P:\alpha(x)=\alpha\} =P\cap f^{-n}P\cap B^s_{n,\alpha},$$
where $B^s_{n,\alpha}=B^s_{n,x}$  is the component of $B\cap
f^{-n}B$ containing $x$. The sets  $T(\alpha)$ actually coincide
with  the  equivalence  classes introduced above.
 
In what follows we assume that $\mu$ is the measure of maximal
entropy of the polynomial diffeomorphism $f$.  Since $P$ has the
product structure, it is homeomorphic to  $P^s\times P^u$, were
$P^s $ and $P^u$ are cross sections.  It was proved in [BLS] that
$\mu$ is a  product measure with respect  to this topological
structure, i.e.\  $\mu\contract P=\mu^-|_{P^s}\otimes\mu^+|_{P^u}$.
 
\proclaim Lemma 5.  $\mu(T(\alpha))\le\mu(P)d^{-n}$.
 
\give Proof.  By Lemma 4, $T$ is naturally homeomorphic to
$T^s\times T^u$, where  $T^s=T^s(x)$ and $T^u=T^u(x)$ be the cross
sections
  of $T=T(\alpha)$ through the point $x\in T(\alpha)$ .
    By the product property of $\mu$,
$$\mu(T)=\mu^-(T^s)\mu^+(T^u).\eqno(10)$$ Since $T^s\subset
P^s(x)$,  $$\mu^-(T^s)\le\mu^-(P^s(x))=\mu^-(P^s).\eqno(11)$$ By
(9), $f^nT^u\subset W^u_{loc}(f^nx)\cap P=P^u(f^nx).$ By the
transformation rule for $\mu^+$ we conclude that $$\mu^+(T^u)\le
d^{-n}\mu^+(P^u).\eqno(12)$$ Now the result follows from
(10)--(12).\qed
 
In the following proof we let ${\rm SFix}_n$ be the set of saddle
periodic points of period dividing $n$.
 
\proclaim Lemma 6.  For any $\epsilon>0$, there exists $C>0$
depending on the Pesin box $P$ such that
$$\liminf_{n\to\infty}{1\over d^{n}}\#\{\alpha\in  {\rm SPer}_n :
{\rm dist}(\alpha,P)<C e^{-n\theta}\}\ge\mu(P).$$
 
\give Proof.  In our work above, we have assigned to any returning
point  $x\in P\cap f^{-n}P$ a periodic point $\alpha(x)\in {\rm
SFix}_n$ exponentially close to $P$. Let $A_n$ denote the set of
periodic points obtained in this way.   Hence, by the Lemma 5,
$$\mu(P) d^{-n}\#A_n\ge\sum_{\alpha\in A_n}\mu(T(\alpha))=
        \mu(P\cap f^{-n}P) $$ and by the mixing property of $\mu$
[BS],  $${\#A_n\over d^n}\ge {\mu(P\cap f^{-n}P)\over
\mu(P)}\to\mu(P)$$ as $n\to\infty$.

Thus we have shown that  $$\liminf_{n\to\infty}{1\over
d^{n}}\#\{\alpha\in {\rm SFix}_n: {\rm dist}(\alpha,P)<C
e^{-n\theta}\}\ge\mu(P).$$ It remains to show that ${\rm SFix}_n$
can be replaced by ${\rm SPer}_n$.  For this, we note that
$\{k:k|n,k<n\}\subset\{k\le n/2\}$.  Since $\#{\rm Fix}_n\le d^n$,
it follows that  $$\#{\rm Fix}_n-\#{\rm Per}_n\le {n\over
2}d^{n/2}.$$ Since this is $o(d^n)$, we see that almost all points
of ${\rm SFix}_n$ have period precisely $n$. \qed

\give Proof of Theorem 1.  Let $\nu_n=d^{-n}\sum_{a\in {\rm
Per}_n}\delta_a$. Consider any limit measure $\nu=\lim\nu_{n(i)}$
of a subsequence of these measures.  It follows from Lemma 6 that
for any Pesin box, $\nu(P)\ge\mu(P)$.  Now let $G$ be any open
set.  Consider a compact set  $K\subset G$ such that
$\mu(G-K)<\epsilon$.
 Let $\eta={\rm dist}(K,\partial G)$. By  Lemma 1 we can cover all
but a set of measure $\epsilon$ by
 disjoint Pesin boxes $P_i$ of size less than $\eta$.  Let $I$
denote the set of those Pesin boxes which intersect  $K$.   Then
$$\nu(G)\ge\sum_{i\in I}\nu(P_i)\ge\sum_{i\in I}\mu(P_i)
\ge\mu(K)-\epsilon\ge\mu(G)-2\epsilon.$$ Hence $\nu\ge \mu$, and
as the both measures are normalized, $\nu=\mu$. \qed

Now  we turn to the Lyapunov exponent $\Lambda(f)$.  For $x\in \cR$ we have
$$\lim_{n\to\infty}{1\over n}\log||Df^n(x)||
=\lim_{n\to\infty}{1\over n}\log\left|Df(x)|_{E^u_x}\right|
=\lim_{n\to\infty}{1\over n}\sum_{j=0}^{n-1}
\log\left| Df(f^jx)|_{E^u_{f^jx}}\right|.$$
Thus with the notation
$$\psi(x)=\log\left|Df(x)|_{E^u_x}\right|\eqno(13)$$
we have
$$\Lambda(f)=\int\psi(x)\mu(x).\eqno(14)$$
 
\proclaim Lemma 7. Let $P$ denote a Pesin box, and let $\epsilon>0$
be given.  Then there exists $N$ sufficiently large that the
periodic saddle point $\alpha(x)$ generated by a returning point
$x\in P\cap f^{-n}P$ with $n\ge N$ satisfies
$dist(E^u_{\alpha(x)},E^u(x))<\epsilon$.
 
\give Proof.  By the construction above, the unstable manifold
$W^u(\alpha(x))$ goes across the topological bidisk $B$, and
$W^u(\alpha(x))\subset B^u_n$.  Since the stable cross section of $B^u_n$
is exponentially small in $n$, we have the desired estimate on the
distance of the tangent spaces. \qed
 
\give Proof of Theorem 2.  The measures $d^{-n}\sum_{a\in
P_n}\delta_a$ converge to $\mu$.  To evaluate the integral in (14), we
may work on a countable, disjoint family of Pesin boxes $P_j$, with the property
that $Df$ and the tangent spaces $E^u$, and thus the expression $\psi$ in (13),
vary by no more than $\epsilon$ on each $P_j$.  By Lemma 7, as $n\to\infty$, the
tangent spaces $E^u_{\alpha(a)}$ converge to within $\epsilon$ of the tangent
spaces $E^u$ given by the Oseledec Theorem.  Thus $\psi(\alpha(x))$ is within
$\epsilon$ of the value of $\psi$ on $P_j$.  It follows that the expression in
Theorem 2 will be within $\epsilon$ of the integral (14). \qed
 
 
\proclaim Corollary 3.  If $P_n$ denotes any of the sets defined in
(1), then $$\Lambda(f)=\lim_{n\to\infty}d^{-n}\sum_{p\in
P_n}\psi(p).$$

\give Proof.  Since  $K$ is compact, it follows that $||Df_x||\le
M$ for all $x\in K$.  And since $P_n\subset K$, it follows that
that $\psi(p)\le M$. Further, the Jacobian  of $f$ is a
constant $a$, so even at a sink orbit of order $n$, we have
$\psi(p)\ge \log\sqrt{|a|}$.  Thus we may use any of the three sets
$P_n$ in (1) in defining the limit.
 \qed
 
\give Proof of Theorem 3. The conclusion of the theorem is a
consequence of the following statement: For any neighborhood $U$
of $c$ there is a number $N_U$ so that for each $n\ge N_U$ there
is a $c\in U$ so that $f_c$ has a sink of period $n$. We will show
that if the previous assertion does not hold then $\Lambda(f_c)$
is harmonic. Thus assume that there exists a neighborhood $U$ and
an infinite sequence $n_i$ such that $f_c$ has no sink of period
$n_i$ for any $c\in U$ and any $i$.
 
Fix an $n=n_i$. Let $V=\{(c,p)\in U\times{\C^2}:f^n_c(p)=p\}$. The
set $V$ is a one dimensional analytic variety with a projection
onto $U$. We can remove from $V$ a discrete set of points
$\{(c_1,p_1),(c_2,p_2),\dots\}$  corresponding to ``bifurcations''
where either $V$ is singular or the projection onto the first
coordinate has a singularity.  Let
$V'=V-\{(c_1,p_1),(c_2,p_2),\dots\}$ and let
$U'=U-\{c_1,c_2,\dots\}$. On $V'$ the period of a periodic point
is constant on each component. Remove components for which the
period is less than $n$. Call the resulting set $V''$.
 
Assume first that we are in the dissipative case, and there are no
sinks of period $n$.  Since there is no sink, the modulus of the larger
eigenvalue must be at least as large as 1, and since it is dissipative, the
modulus of the smaller eigenvalue must be no larger than $\sqrt|a| < 1$.  Thus
at each point $(c,p)$ in $V'$  there is a unique largest eigenvalue for
$Df_c^n(p)$ call it $\lambda^+(c,p)$.  The function $\lambda^+$ is a
continuously chosen root of the characteristic equation of $Df_c^n(p)$ so it is
holomorphic. In particular the function $\log|\lambda^+(c,p)|$ is harmonic on
$V'$ and the function $$U'\ni c\mapsto {1\over d^n}\sum_{p\in P_n}{1\over
n}\log|\lambda^+(c,p)|={1\over d^n}\sum_{p\in P_n}\psi(p)$$  is harmonic on
$U'$. This function extends to a harmonic function $\Lambda_n$ on $U$.  Now
Theorem 2 implies that $\Lambda_n$ converges pointwise to $\Lambda$. Since each
function $\Lambda_n$ is harmonic we conclude that $\Lambda$ is harmonic. This
proves our assertion in the dissipative case.
 
Assume now that we are in the volume preserving case. Note that
since the function $c\mapsto \det f_c$ is holomorphic and of
constant norm it is actually constant. Assume that there are no
Siegel balls of period $n$.  By the two-dimensional version of the Siegel
linearization theorem (see [Z]),
this implies that there are no elliptic periodic
points of period $n$ for which the eigenvalues satisfy certain Diophantine
conditions.  In particular if at any elliptic periodic point the eigenvalues
are not constant as a function of $c$ then we can vary the parameter so that
$\lambda_1$ varies through an interval on the unit circle. This implies that
there is some parameter value at which the Diophantine conditions are
verified for the eigenvalues $\lambda_1$ and
$\lambda_2=const/\lambda_1$. We conclude that at any elliptic
point $(c,p)$ of period $n$ the eigenvalues are locally
independent of $c$. This implies that the eigenvalues are constant
on the component of $V'$ which contains the point $(c,p)$.

As before we remove from $V'$ the components on which the period
is less than $n$. In addition we remove components on which both
eigenvalues are constant and have modulus 1. The remaining variety
$V''$ consists of
 saddles. Arguing as before we see that $\Lambda_n$ is harmonic
hence $\Lambda$ is harmonic. This completes the proof of the
theorem.  \qed
 
\give Proof of Theorem 4. Let $J^*$ denote the support of $\mu$, which is the
closure of the saddle points (cf.\ [BS3]). If the maps $f_c|_{J^*}$ are
topologically conjugate, then the conjugacy preserves the set of periodic
points, so that each periodic point $p$  is part of a selection $c\mapsto
p(c)$.  As in the proof of Theorem 3, each function $\Lambda_n$ is harmonic;
hence $\Lambda$ is harmonic. \qed

\centerline{\bf References}
 
\item{[BLS]} E. Bedford, M. Lyubich, and J. Smillie, Polynomial diffeomorphisms
of $\cx2$. IV: The measure of maximal entropy and laminar currents. Invent.
Math., to appear.
 
\item{[BS1]}  E. Bedford and J. Smillie, Polynomial diffeomorphisms of $\cx2$:
Currents, equilibrium measure and hyperbolicity. Invent. Math. 87, 69--99
(1990)
 
\item{[BS2]}  E. Bedford and J. Smillie, Polynomial diffeomorphisms of $\cx2$
II: Stable manifolds and recurrence. J. AMS 4, 657--679 (1991)
 
\item{[BS3]}  E. Bedford and J. Smillie, Polynomial diffeomorphisms of $\cx2$
III: Ergodicity, exponents and entropy of the equilibrium measure. Math. Ann.
294. 395--420 (1992)
 
\item{[B1]}  R. Bowen, Periodic points and measures for axiom A diffeomorphisms.
 Trans. AMS. 154, 377-397 (1971).
 
\item{[B2]}  R. Bowen, Equilibrium states and the ergodic theory of Anosov
diffeomorphisms, (Lect. Notes, Math., vol 470) Berlin Heidelberg New York:
Springer 1975
 
\item{[Br]}  H. Brolin, Invariant sets under iteration of rational functions.
Ark. Mat. 6, 103--144 (1965)
 
\item{[FS]}  J.-E. Forn\ae ss and N. Sibony, Complex H\'enon mappings in $\cx2$
and Fatou Bieberbach domains. Duke Math. J. 65, 345--380 (1992)
 
\item{[K]}  A. Katok, Lyapunov exponents, entropy and periodic orbits for
diffeomorphisms.  Publ. Math. Inst. Hautes Etud. Sci. 51, 137--174 (1980)
 
\item{[L]}  M. Lyubich, Entropy of analytic endomorphisms of the Riemannian
sphere. Funct. Anal. Appl. 15, 300--302 (1981)
 
\item{[PS]}  C. Pugh and M. Shub, Ergodic attractors. Trans. AMS 312, 1--54
(1989)
 
\item{[Si]}  N. Sibony, Iteration of polynomials, U.C.L.A. course lecture notes.

\item{[S]}  J. Smillie, the entropy of polynomial diffeomorphisms of $\cx2$.
Ergodic Theory Dyn. Syst. 10, 923--827 (1990)
 
\item{[T]} P. Tortrat, Aspects potentialistes de l'it\'eration des
polyn\^omes.  In: S\'eminaire de Th\'eorie du Potentiel Paris, No. 8 (Lect.
Notes Math., vol. 1235) Berlin Heidelberg New York: Springer 1987.
 
\item{[Z]}  E. Zehnder, A simple proof of a generalization of a Theorem by C.
L. Siegel. Geometry and Topology (J. Palis and M. do Carmo, eds.), Lecture
Notes in Math., vol. 597, Springer-Verlag, New York, 855--866 (1977)

\bigskip \noindent Indiana University, Bloomington, IN 47405
\medskip \noindent SUNY at Stony Brook, Stony Brook, NY 11794
\medskip \noindent Cornell University, Ithaca, NY 14853

 \bye